\documentclass{amsart}
\usepackage{amsfonts}
\usepackage{amssymb}
\usepackage{amsmath}
\usepackage{amsthm}

\newtheorem{theorem}{Theorem}

\newtheorem{lemma}{Lemma}

\theoremstyle{definition}
\newtheorem{definition}{Definition}

\begin{document}

\title{Holomorphic Motions and Normal Forms in
Complex Analysis}

\author
{Yunping Jiang}
\address{Department of Mathematics\\
Queens College of the City University of New York\\
Flushing, NY 11367-1597\\
and\\
Department of Mathematics\\
Graduate School of the City University of New York\\
365 Fifth Avenue, New York, NY 10016\\
and\\
Academy of Mathematics and System Sciences\\
Chinese Academy of Sciences, Beijing 100080}
\email[Jiang]{yungc@forbin.qc.edu}



\thanks{The research is partially supported by
NSF grants and PSC-CUNY awards and the Hundred Talents Program
from Academia Sinica and to be appeared in the Proceedings
of ICCM2004, International Press}

\begin{abstract}
We give a brief review of holomorphic motions and its relation
with quasiconformal mapping theory. Furthermore, we apply the
holomorphic motions to give new proofs of famous K\"onig's Theorem
and B\"ottcher's Theorem in classical complex analysis.
\end{abstract}

\maketitle

\section{Introduction}

Two of fundamental tools in the study of complex dynamical systems
are K\"onig's Theorem and B\"otthcher's Theorem in classical
complex analysis, which were proved back to 1884~\cite{Kon} and
1904~\cite{Bot}, respectively, by using some well-known methods in
complex analysis. These theorems say that an attractive or
repelling or super-attractive fixed point of an analytic map can
be written into a normal form under suitable conformal changes of
coordinate.

During the study of complex dynamical systems, a subject called
holomorphic motions becomes more and more interesting and useful.
The subject of holomorphic motions over the open unit disk shows
some interesting connections between classical complex analysis
and problems on moduli. This subject even becomes an interesting
branch in complex analysis~\cite{AM,BR,Do,EM,Lie,Mit,ST,Sl}.

In this paper, we apply the holomorphic motions to show new proofs
of K\"onig's Theorem and B\"otthcher's Theorem. From the technical
point of views, our proofs are more complicate and use a
sophistical result in holomorphic motions. But from the conceptual
point of views, our proofs give some inside mechanism of the
normal forms for fixed points. The technique involved in the new
proofs will develop some new results too (refer to~\cite{Ji}).

In our new proofs, we use quasiconformal theory, in particular,
holomorphic motions over the open unit disk. Holomorphic motions
over a simply connected complex Banach manifold are also
interesting and useful. In sequel, we will discuss them and show
some application to complex dynamical systems~\cite{JM}.

\vskip10pt \noindent {\bf Acknowledgment.} This article was first
discussed in the complex analysis and dynamical systems seminar at
the CUNY Graduate Center in the Spring Semester of 2003. I would
like to thank all participants for their patience. During the
discussion Professors Linda Keen, Fred Gardiner, and Nikola Lakic
provided many useful comments to improve and polish this article.
Sudeb Mitra explained to me several points in the development of
the measurable Riemann mapping theorem and holomorphic motions. I
express my sincere thanks to them.

\section{Holomorphic Motions and Quasiconformal Maps}

In the study of complex analysis, the measurable Riemann mapping
theorem plays an important role. Consider the Riemann sphere
$\hat{\mathbb C}$. A measurable function $\mu$ on $\hat{\mathbb
C}$ is called a Beltrami coefficient if there is a constant $0\leq
k<1$ such that $\|\mu\|_{\infty} \leq k$, where
$\|\cdot\|_{\infty}$ means the $L^{\infty}$-norm of $\mu$ on
$\hat{\mathbb C}$. The equation
$$
H_{\overline{z}}=\mu H_{z}
$$
is called the Beltrami equation with the given Beltrami
coefficient $\mu$. The measurable Riemann mapping theorem says
that the Beltrami equation has a solution $H$ which is a
quasiconformal homeomorphism of $\hat{\mathbb C}$ whose
quasiconformal dilatation is less than or equal to
$K=(1+k)/(1-k)$. The study of the measurable Riemann mapping
theorem has a long history since Gauss considered in 1820's the
connection with the problem of finding isothermal coordinates for
a given surface. As early as 1938, Morrey~\cite{Mo} systematically
studied homeomorphic $L^{2}$-solutions of the Beltrami equation.
But it took almost twenty years until in 1957 Bers~\cite{Be}
observed that these solutions are quasiconformal (refer
to~\cite[pp. 24]{Le}). Finally the existence of a solution to the
Beltrami equation under the most general possible circumstance,
namely, for measurable $\mu$ with $\|\mu\|_{\infty}<1$, was shown
by Bojarski~\cite{Bo}. In this generality the existence theorem is
sometimes called the measurable Riemann mapping theorem (refer
to~\cite[pp. 10]{GL}).

If one only considers a normalized solution in the Beltrami
equation (a solution fixes $0$, $1$, and $\infty$), then $H$ is
unique, which is denoted as $H^{\mu}$. The solution $H^{\mu}$ is
expressed as a power series made up of compositions of singular
integral operators applied to the Beltrami equation on the Riemann
sphere. In this expression, if one considers $\mu$ as a variable,
then the solution $H^{\mu}$ depends on $\mu$ analytically. This
analytic dependence was emphasized by Ahlfors and Bers in their
1960 paper~\cite{AB} and is essential in determining a complex
structure for Teichm\"uller space (refer
to~\cite{Al,GL,Le,Li,Na}). Note that when $\mu\equiv 0$, $H^{0}$
is the identity map. A $1$-quasi-conformal map is conformal.
Twenty years later, due to the development of complex dynamics,
this analytic dependence presents an even more interesting
phenomenon called holomorphic motions as follows.

Let $\Delta_{r}=\{ c\in {\mathbb C} \;|\; |c|<r\}$ be the disk
centered at $0$ and of radius $r>0$. In particular, we use
$\Delta$ to denote the unit disk. Given a Beltrami coefficient
$\mu$, consider a family of Beltrami coefficients $c\mu$ for $c\in
\Delta$ and the family of normalized solutions $H^{c\mu}$. Note
that $H^{c\mu}$ is a quasiconformal homeomorphism whose
quasiconformal dilatation is less than or equal to
$(1+|c|k)/(1-|c|k)$. Moreover, $H^{c\mu}$ is a family which is
holomorphic on $c$. Consider a subset $E$ of $\hat{\mathbb C}$ and
its image $E_{c}=H^{c\mu}(E)$. One can see that $E_{c}$ moves
holomorphically in $\hat{\mathbb C}$ when $c$ moves in $\Delta$.
That is, for any point $z\in E$, $z(c)=H^{c\mu}(z)$ traces a
holomorphic path starting from $z$ as $c$ moves in the unit disk.
Although $E$ may start out as smooth as a circle and although the
points of $E$ move holomorphically, $E_{c}$ can be an interesting
fractal with fractional Hausdorff dimension for every $c\neq 0$
(see~\cite{GK}).

Surprisingly, the converse of the above fact is true too. This
starts from the famous $\lambda$-lemma of Ma\~n\'e, Sad, and
Sullivan~\cite{MSS} in complex dynamical systems. Let us start to
understand this fact by first defining holomorphic motions.

\medskip
\begin{definition}[Holomorphic Motions]~\label{hm} Let $E$ be a
subset of $\hat{\mathbb C}$. Let
$$
h (c, z) : \Delta_{r}\times E\to \hat{\mathbb C}
$$
be a map. Then $h$ is called a holomorphic motion of $E$
parametrized by $\Delta_{r}$ if
\begin{enumerate}
\item $h (0, z)=z$ for $z\in E$;
\item for any fixed $c\in
\Delta_{r}$, $h (c, \cdot): E\to \hat{\mathbb C}$ is injective;
\item
for any fixed $z$, $h (\cdot,z): \Delta_{r} \to \hat{\mathbb C}$
is holomorphic.
\end{enumerate}
\end{definition}

\medskip
For example, for a given Beltrami coefficient $\mu$,
$$
H(c, z)=H^{c\mu}(z): \Delta\times \hat{\mathbb C}\to \hat{\mathbb
C}
$$
is a holomorphic motion of $\hat{\mathbb C}$ parametrized by
$\Delta$.

Note that even continuity does not directly enter into the
definition; the only restriction is in the $c$ direction. However,
continuity is a consequence of the hypotheses from the proof of
the $\lambda$-lemma of Ma\~n\'e, Sad, and Sullivan~\cite[Theorem
2]{MSS}. Moreover, Ma\~n\'e, Sad, and Sullivan prove in~\cite{MSS}
that

\medskip
\begin{lemma}[$\lambda$-Lemma]~\label{ll}
A holomorphic motion of a set $E\subset \hat{\mathbb C}$
parametrized by $\Delta_{r}$ can be extended to a holomorphic
motion of the closure of $E$ parametrized by the same
$\Delta_{r}$.
\end{lemma}

\medskip
Furthermore, Ma\~n\'e, Sad, and Sullivan show in~\cite{MSS} that
$f(c, \cdot)$ satisfies the Pesin property. In particular, when
the closure of $E$ is a domain, this property can be described as
the quasiconformal property. A further study of this
quasiconformal property is given by Sullivan and
Thurston~\cite{ST} and Bers and Royden~\cite{BR}. In~\cite{ST},
Sullivan and Thurston prove that there is a universal constant
$a>0$ such that any holomorphic motion of any set $E\subset
\hat{\mathbb C}$ parametrized by the open unit disk $\Delta$ can
be extended to a holomorphic motion of $\hat{\mathbb C}$
parametrized by $\Delta_{a}$. In~\cite{BR}, Bers and Royden show,
by using classical Teichm\"uller theory, that this constant
actually can be taken to be $1/3$. Moreover, in the same paper,
Bers and Royden show that in any holomorphic motion $H(c,z):
\Delta_{r}\times \hat{\mathbb C}\to \hat{\mathbb C}$ for $0<r\leq
1$, $H(c,\cdot): \hat{\mathbb C}\to \hat{\mathbb C}$ is a
quasiconformal homeomorphism whose quasiconformal dilatation less
than or equal to $(1+|c|)/(1-|c|)$ for $c\in \Delta_{r}$. In the
both papers~\cite{ST,BR}, they expect $a=1$. This was eventually
proved by Slodkowski in~\cite{Sl}.

\medskip
\begin{theorem}[Slodkowski's Theorem]~\label{st}
Suppose
$$
h(c,z): \Delta\times E\to \hat{\mathbb C}
$$
is a holomorphic motion of a set $E\subset \hat{\mathbb C}$
parametrized by $\Delta$ . Then $h$ can be extended to a
holomorphic motion
$$
H(c, z): \Delta\times \hat{\mathbb C}\to \hat{\mathbb C}
$$
of ${\mathbb C}$ parametrized by also $\Delta$. Moreover,
following~\cite[Theorem 1]{BR}, for every $c\in \Delta$, $H(c,
\cdot): \hat{\mathbb C}\to \hat{\mathbb C}$ is a quasiconformal
homeomorphism whose quasiconformal dilatation
$$
K(H(c,\cdot)) \leq \frac{1+|c|}{1-|c|}.
$$
\end{theorem}

\medskip
Holomorphic motions of a set $E\subset \hat{\mathbb C}$
parametrized by a connected complex manifold with a base point can
be also defined. They have many interesting relationships with the
Teichm\"uller space $T(E)$ of a closed set $E$ (refer
to~\cite{Mit}).

In addition to the references we mentioned above, there is a
partial list of references~\cite{As,DHOP,EKK,MS,Po} about
holomorphic motions and quasiconformal mapping theory. The reader
who is interested in holomorphic motions may refer to those papers
and books.

\section{A new proof of K\"onig's Theorem}

We first give a new proof of K\"onig's Theorem. The idea of the
new proof follows the viewpoint of holomorphic motions. For the
classical proof of K\"onig's Theorem, the reader may refer to
K\"onig's original paper~\cite{Kon} or most recent
books~\cite{CG,Mi}. As we mentioned in the introduction, from the
technique point of views, the new proof  may be more complicate,
but from the conceptual point of views, it gives some inside
mechanism of the normal form for an attractive fixed point.

\medskip
\begin{theorem}[K\"onig's Theorem]~\label{kth}
Let $f(z)=\lambda z+\sum_{j=2}^{\infty} a_{j}z^{j}$ be an analytic
function defined on $\Delta_{r_{0}}$, $r_{0}>0$. Suppose
$0<|\lambda|<1$ or $|\lambda|>1$. Then there is a conformal map
$\phi: \Delta_{\delta} \to \phi(\Delta_{\delta})$ for some
$0<\delta <r_{0}$ such that
$$
\phi^{-1}\circ f\circ \phi (z) =\lambda z.
$$
The conjugacy $\phi^{-1}$ is unique up to multiplication of
constants.
\end{theorem}

\medskip
\begin{proof}
We only need to prove it for $1<|\lambda|<1$. In the case of
$|\lambda|>1$, we can consider $f^{-1}$.

First, we can find a $0<\delta<r_{0}$ such that
$$
|f(z)|< |z|, \quad z\in \overline{\Delta}_{\delta}
$$
and $f$ is injective on $\overline{\Delta}_{\delta}$. For every
$0<r\leq \delta$, let
$$
S_{r}=\{ z\in {\mathbb C}\;|\; |z|=r\}
$$ \
and
$$
T_{r}=|\lambda| S_{r}=\{ z\in {\mathbb C}\;|\; |z|=|\lambda| r\}.
$$
Denote $E=S_{r}\cup T_{r}$. Define
$$
\phi_{r} (z) =\left\{
\begin{array}{ll}
        z, & z\in S_{r}; \cr
        f(\frac{z}{\lambda}), & z\in T_{r}.
\end{array}\right.
$$
It is clear that
$$
\phi_{r}^{-1}\circ f\circ\phi_{r} (z)=\lambda z
$$
for $z\in S_{r}$.

Now write $\phi_{r}(z)=z\psi_{r} (z)$ for $z\in T_{r}$, where
$$
\psi_{r} (z) =1+\sum_{j=1}^{\infty} \frac{a_{j+1}}{\lambda^{j+1}}
z^{j}. $$ Define
$$
h_{r}(c, z) =\left\{
\begin{array}{ll}
z, & z\in S_{r}\cr z\psi_{r} (\frac{\delta cz}{r}), & z\in T_{r}
\end{array}\right.
: \Delta\times E\to \hat{\mathbb C}.
$$
Note that
$$
h(c,z) = z \psi_{r} \Big( \frac{cz\delta}{r}\Big)=\frac{r}{c\delta
}\phi \Big( \frac{cz\delta}{r}\Big) =\frac{r}{c\delta} f\Big(
\frac{cz\delta}{r\lambda}\Big), \quad z\in T_{r}, c\neq 0.
$$
For each fixed $z\in E$, it is clear that $h (c, z)$ is a
holomorphic function of $c\in \Delta$. For each fixed
$c\in\Delta$, the restriction $h (c, \cdot)$ to $S_{r}$ and
$T_{r}$, respectively, are injective. Now we claim that their
images do not cross either. That is because for any $z\in T_{r}$,
$|z|=|\lambda|r$ and $|cz\delta|/|r\lambda| \leq \delta$, so
$$
|h(c, z)|=\Big| \frac{r}{c\delta}\Big| \Big| f\Big( \frac{
cz\delta}{r \lambda}\Big) \Big| <\Big| \frac{r}{c\delta}\Big|
\Big| \frac{cz\delta}{r\lambda}\Big| =r.
$$
Therefore, $h (c,z): \Delta\times E\to \hat{\mathbb C}$ is a
holomorphic motion because we also have $h(0, z) =z$ for all $z\in
E$. From Slodkowski's Theorem, $h$ can be extended to a
holomorphic motion $H(c,z): \Delta\times \hat{\mathbb C}\to
\hat{\mathbb C}$, and moreover, for each fixed $c\in \Delta$,
$H_{c}=h(c,\cdot):\hat{\mathbb C}\to \hat{\mathbb C}$ is a
quasiconformal homeomorphism whose quasiconformal dilatation is
less than or equal to $(1+|c|)/(1-|c|)$. Now take $c_{r}=r/\delta$
and consider $H(c_{r},\cdot)$. We have $H(c_{r},
\cdot)|E=\phi_{r}$. Let
$$
A_{r,j}= \{z\in {\mathbb C}\;|\; |\lambda|^{j+1}r\leq |z|\leq
|\lambda|^{j} r\}.
$$
We still use $\phi_{r}$ to denote $H(c_{r},\cdot)|A_{r, 0}$.

For an integer $k>0$, take $r=r_{k}=\delta |\lambda|^{k}$. Then
$$
\overline{\Delta}_{\delta} =\cup_{j=-k}^{\infty} A_{r,j}\cup
\{0\}.
$$
Extend $\phi_{r}$ to $\overline{\Delta}_{\delta}$, which we still
denote as $\phi_{r}$, as follows.
$$
\phi_{r}(z) = f^{-j}(\phi_{r}((\lambda^{n}z)), \quad z\in A_{r,j},
\quad j=-k, \cdots, -1, 0, 1, \cdots,
$$
and $\phi_{r}(0)=0$. Since $\phi_{r}|E$ is a conjugacy from $f$ to
$\lambda z$, $\phi_{r}$ is continuous on $\Delta_{\delta}$. Since
$f$ is conformal, $\phi_{r}$ is quasiconformal whose
quasiconformal dilatation is the same as that of $H(c_{r}, \cdot)$
on $A_{r,0}$. So the quasiconformal dilatation of $\phi_{r}$ on
$\Delta_{\delta}$ is less than or equal to $(1+r)/(1-r)$.
Furthermore,
$$
f (\phi_{r} (z)) =\phi_{r}(\lambda z), \quad z\in \Delta_{\delta}.
$$

Since $f(z) =\lambda z( 1+O(z))$, $f^{k}(z) =\lambda^{k} z
\prod_{i=0}^{k-1} (1+O(\lambda^{i}z))$. Because
$|\lambda|^{k}r_{k}=\delta$, the range of $\phi_{r_{k}}$ on
$\Delta_{\delta}$ is a Jordan domain bounded above and below
uniformly on $k$. In addition, $0$ is fixed by $\phi_{k}$ and the
quasiconformal dilatations of the $\phi_{k}$ are uniformly
bounded. Therefore, the sequence $\{
\phi_{r_{k}}\}_{k=1}^{\infty}$ is a compact family
(see~\cite{Al}). Let $\phi$ be a limiting map of this family. Then
we have
$$
f(\phi (z)) =\phi(\lambda z), \quad z\in \Delta_{\delta}.
$$
The quasiconformal dilatation of $\phi$ is less than or equal to
$(1+r_{k})/(1-r_{k})$ for all $k>0$. So $\phi$ is a
$1$-quasiconformal map, and thus is conformal. This is the proof
of the existence.

For the sake of completeness, we also provide the proof of
uniqueness but this is not new and the reader can find it
on~\cite{CG,Mi}. Suppose $\phi_{1}$ and $\phi_{2}$ are two
conjugacies such that
$$
\phi^{-1}_{1}\circ f\circ \phi_{1} (z) =\lambda z \quad
\hbox{and}\quad \phi^{-1}_{2}\circ f\circ \phi_{2} (z) =\lambda z,
\quad z\in \Delta_{\delta}.
$$
Then for $\Phi=\phi_{2}^{-1}\circ \phi_{1}$, we have $\Phi(\lambda
z)=\lambda \Phi(z)$. This implies that $\Phi'(\lambda z)=\Phi'(z)$
for any $z\in \Delta_{\delta}$. Thus
$\Phi'(z)=\Phi'(\lambda^{n}z)=\Phi(0)$. So $\Phi(z)=\hbox{const}$
and $\phi_{2}^{-1}=\hbox{const.}\cdot \phi_{1}^{-1}$.
\end{proof}

\medskip
\section{A new proof of B\"ottcher's Theorem}

In this section, we give a new proof of B\"ottcher's Theorem. The
idea of the new proof follows the viewpoint of holomorphic
motions. For the classical proof of B\"ottcher's Theorem, the
reader may refer to B\"ottcher's original paper~\cite{Bot} or most
recent books~\cite{CG,Mi}. The idea of the proof is basically the
same as that in the previous section, but the actual proof is
little bit different. The reason is that in the previous case, $f$
is a homeomorphism so we can iterate both forward and backward,
but in B\"ottcher's Theorem, $g$ is not a homeomorphism.

\medskip
\begin{theorem}[B\"ottcher's Theorem] Suppose $f(z)
=\sum_{j=n}^{\infty}a_{j}z^{j}$, $a_{n}\neq 0$, $n\geq 2$, is
analytic on a disk $\Delta_{\delta_{0}}$, $\delta_{0}>0$. Then
there exists a conformal map $\phi: \Delta_{\delta}\to
\phi(\Delta_{\delta})$ for some $\delta>0$ such that
$$
\phi^{-1}\circ f\circ \phi (z) =z^{n}, \quad z\in \Delta_{\delta}.
$$
The conjugacy $\phi^{-1}$ is unique up to multiplication by
$(n-1)^{th}$-roots of the unit.
\end{theorem}

\medskip
\begin{proof}
Conjugating by $z\to bz$, we can assume $a_{n}=1$, i.e.,
$$
f(z) =z^{n} +\sum_{j=n+1}^{\infty} a_{j}z^{j}.
$$

We use $\Delta_{r}^{*}=\Delta_{r}\setminus \{0\}$ to mean a
punctured disk of radius $r>0$. Write
$$
f(z)=z^{n}(1+\sum_{j=1}^{\infty} a_{j+n}z^{j}).
$$
Assume $0<\delta_{1}<\min\{ 1/2, \delta_{0}/2\}$ is small enough
such that
$$
1+\sum_{j=1}^{\infty} a_{j+n}z^{j} \neq 0 \quad \hbox{and} \quad
\frac{1}{\hbox{$\root n \of  {|1+\sum_{j=1}^{\infty}
a_{j+n}z^{j}}|$}} \geq \frac{1}{2}, \quad z\in
\Delta_{2\delta_{1}}.
$$
Then $f: \Delta_{2\delta_{1}}^{*}\to f(\Delta_{2\delta_{1}}^{*})$
is a covering map of degree $n$.

Let $0<\delta < \delta_{1}$ be a fixed number such that
$f^{-1}(\Delta_{\delta}) \subset \Delta_{\delta_{1}}$. Since
$$
z\to z^{n}: \Delta_{\hbox{$\root n \of {\delta}$}}^{*} \to
\Delta_{\delta}^{*} \quad \hbox{and}\quad  f:
f^{-1}(\Delta_{\delta}^{*}) \to \Delta_{\delta}^{*}
$$
are both of covering maps of degree $n$, the identity map of
$\Delta_{\delta}$ can be lifted to a holomorphic diffeomorphism
$$
h: \Delta_{\hbox{$\root n \of {\delta}$}}^{*}\to
f^{-1}(\Delta_{\delta}^{*}),
$$
i.e., $h$ is a map such that the diagram
$$
\begin{array}{ccc}
\Delta_{\hbox{$\root n \of {\delta}$}}^{*} & {\buildrel h \over
\longrightarrow}& f^{-1}(\Delta_{\delta}^{*})\cr \downarrow z\to
z^{n} &         &\downarrow f\cr
\Delta_{\delta}^{*} &{\buildrel
\hbox{id} \over \longrightarrow}&
   \Delta_{\delta}^{*}
\end{array}
$$
commutes. We pick the lift so that
$$
h (z) =z \Big( 1+\sum_{j=2}^{\infty} b_{j}z^{j-1}\Big) =z\psi(z).
$$
From
$$
f(h(z))=z^{n}, \quad z\in \Delta_{\hbox{$\root n \of
{\delta}$}}^{*},
$$
we get
$$
|h(z)| =\frac{|z|}{\hbox{$\root n \of {|1+\sum_{j=1}^{\infty}
a_{n+j}(h(z))^{j}|}$}} \geq \frac{|z|}{2}.
$$

For any
$$
0<r\leq \max \Big\{ \Big( \frac{1}{2}\Big)^{\frac{n}{(n-1)}},
\delta^{n}\Big\},
$$
let $S_{r}=\{ z\in {\mathbb C}\; |\; |z|=r\}$ and $T_{r}=\{ z\in
{\mathbb C}\;|\; |z|=\hbox{$\root n \of r$}\}$. Consider the set
$E=S_{r}\cup T_{r}$ and the map
$$
\phi_{r}(z) =\left\{
\begin{array}{ll}
z, & z\in S_{r}\cr
z\psi(z), & z\in T_{r}.
\end{array}\right.
$$
Define
$$
h_{r}(c, z) =\left\{
\begin{array}{ll}
z, & z\in S_{r}\cr z\psi \Big( \frac{cz}{\hbox{$\root n \of
r$}}\Big), & z\in T_{r}
\end{array}\right.
:\Delta\times E\to \hat{\mathbb C}.
$$
Note that
$$
z\psi\Big( \frac{cz}{\hbox{$\root n \of r$}}\Big)
=\frac{\hbox{$\root n \of r$}}{c} h \Big( \frac{cz}{\hbox{$\root n
\of r$}}\Big), \quad z\in T_{r},\; c\not= 0.
$$
This implies that
$$
|h_{r}(c,z)| =\frac{\hbox{$\root n \of r$}}{|c|} h\Big(
\frac{cz}{\hbox{$\root n \of r$}}\Big) \Big| \geq
\frac{\hbox{$\root n \of r$}}{|c|} \frac{|cz|}{2 \hbox{$\root n
\of r$}} \geq \frac{\hbox{$\root n \of r$}}{2}>r, \quad z\in
T_{r}.
$$
So images of $S_{r}$ and $T_{r}$ under $h_{r}(c,z)$ do not cross
each other.

Now let us check that $h_{r}(c,z)$ is a holomorphic motion. First
$h_{r}(0,z)=z$ for $z\in E$. For fixed $x\in E$, $h_{r}(c, z)$ is
holomorphic on $c\in \Delta$. For fixed $c\in \Delta$,
$h_{r}(c,z)$ restricted to $S_{r}$ and $T_{r}$, respectively, are
injective. But the images of $S_{r}$ and $T_{r}$ under
$h_{r}(c,z)$ do not cross each other. So $h_{r}(c, z)$ is
injective on $E$.  Thus
$$
h_{r}(c,z): \Delta\times E\to \hat{\mathbb C}
$$
is a holomorphic motion. By Slodkowski's Theorem, it can be
extended to a holomorphic motion
$$
H_{r}(c, z): \Delta\times \hat{\mathbb C}\to \hat{\mathbb C}.
$$
And moreover, for each $c\in \Delta$, $H_{r}(c,\cdot)$ is a
quasiconformal map whose quasiconformal dilatation satisfies
$$
K(H(c, \cdot)) \leq \frac{1+|c|}{1-|c|}.
$$

Now consider $H(\hbox{$\root n \of r$}, \cdot)$. It is a
quasiconformal map with quasiconformal constant
$$
K_{r}\leq \frac{1+\hbox{$\root n \of r$}}{1-\hbox{$\root n \of
r$}}.
$$
Let
$$
A_{r,j}=\{z\in {\mathbb C}\;|\; \hbox{$\root n^{j} \of r$} \leq
|z|\leq \hbox{$\root n^{j+1} \of r$} \},\quad j=0, 1, 2, \cdots.
$$
Consider the restriction $\phi_{r,0}=H(\hbox{$\root n \of r$},
\cdot)|A_{r, 0}$. It is an extension of $\phi_{r}$, i.e.,
$\phi_{r,0}|E=\phi_{r}$.

Let $\tilde{A}_{r,0}$ be the annulus bounded by $S_{r}$ and
$f^{-1}(S_{r})$ and define $\tilde{A}_{r, j} =f^{-j}
(\tilde{A}_{r, 0})$, $j\geq 0$. Since $z\to z^{n}: A_{r,1}\to
A_{r, 0}$ and $f: \tilde{A}_{r,1}\to \tilde{A}_{r,0}$ are both
covering maps of degree $n$, so $\phi_{r,0}$ can be lifted to a
quasiconformal map $\phi_{r, 1}: A_{r, 1}\to \tilde{A}_{r, 1}$,
i.e., the following diagram
$$
\begin{array}{ccc}
A_{r,1} &{\buildrel \phi_{r,1} \over \longrightarrow}&
\tilde{A}_{r,1}\cr \downarrow z\to z^{n} &         &\downarrow
f\cr
A_{r,0} &{\buildrel \phi_{r,0} \over \longrightarrow}&
\tilde{A}_{r,0}
\end{array}
$$
commutes. We pick the lift $\phi_{r, 1}$ such that it agrees with
$\phi_{r,0}$ on $T_{r}$. The quasiconformal dilatation of
$\phi_{r,1}$ is less than or equal to $K_{r}$.

For an integer $k>0$, take $r=r_{k}=\delta^{n^k}$. Inductively, we
can define a sequence of $K_{r}$-quasiconformal maps $\{\phi_{r,
j}\}_{j=0}^{k}$ such that
$$
\begin{array}{ccc}
A_{r,j} &{\buildrel \phi_{r,j} \over \longrightarrow}&
\tilde{A}_{r,j}\cr \downarrow z\to z^{n} &         &\downarrow
f\cr A_{r,j-1} &{\buildrel \phi_{r,j-1} \over \longrightarrow}&
\tilde{A}_{r,j-1}
\end{array}
$$
commutes and $\phi_{r,j}$ and $\phi_{r, j-1}$ agree on the common
boundary of $A_{r,j}$ and $A_{r, j-1}$. Note that
$$
\overline{\Delta}_{\delta}=\Delta_{r}\cup \cup_{j=0}^{k} A_{r,j}.
$$
Now we can define a quasiconformal map, which we still denote by
$\phi_{r}$ as follows.
$$
\phi_{r} (z) =\left\{
\begin{array}{lll} z, & z\in \Delta_{r}; &\cr
\phi_{r, j}, & z\in A_{r, j},& j=0, 1, \cdots, k.
\end{array}\right.
$$
The quasiconformal dilatation of $\phi_{r}$ on $\Delta_{\delta}$
is less than or equal to $K_{r}$ and
$$
f(\phi_{r}(z)) =\phi_{r}(z^{n}), \quad z\in \cup_{j=1}^{k}
A_{r,j}.
$$

Since $f(z) =z^{n}( 1+O(z))$, $f^{k}(z) = z^{n^{k}}
\prod_{i=0}^{k-1} (1+O(z^{n^i}))$. Because ${\root n^{k} \of
r_{k}}=\delta$, the range of $\phi_{r_{k}}$ on $\Delta_{\delta}$
is a Jordan domain bounded above and below uniformly in $k$. In
addition, $0$ is fixed by $\phi_{k}$ and the quasiconformal
dilatations of the $\phi_{k}$ are uniformly bounded in $k$.
Therefore, the sequence $\{ \phi_{r_{k}}\}_{k=1}^{\infty}$ is a
compact family (see~\cite{Al}). Let $\phi$ be a limiting map of
this family. Then we have
$$
f(\phi (z)) =\phi(z^{n}), \quad z\in \Delta_{\delta}.
$$
Since the quasiconformal dilatation of $\phi$ is less than or
equal to $(1+{\root n \of r_{k}}) /(1-{\root n \of r_{k}})$ for
all $k>0$, it follows that $\phi$ is a $1$-quasiconformal map, and
thus conformal. This is the proof of the existence.

Suppose $\phi_{1}$ and $\phi_{2}$ are two conjugacies such that
$$
\phi^{-1}_{1}\circ f\circ \phi_{1} (z) =z^{n} \quad
\hbox{and}\quad \phi^{-1}_{2}\circ f\circ \phi_{2} (z) =z^{n},
\quad z\in \Delta_{\delta}.
$$
For
$$
\Phi (z) =\phi_{2}^{-1}\circ \phi_{1} (z) =\sum_{j=1}^{\infty}
a_{j}z^{j},
$$
we have $\Phi( z^{n})=(\Phi(z))^{n}$. This implies
$a_{1}^{n}=a_{1}$ and $a_{j}=0$ for $j\geq 2$. Since $a_{1}\neq
0$, we have $a_{1}^{n-1}=1$ and
$\phi_{2}^{-1}=a_{1}\phi_{1}^{-1}$. This is the uniqueness.
\end{proof}

\medskip

\bibliographystyle{amsalpha}

\end{document}